%% file: main.tex
\documentclass[a4paper,review]{cas-sc}

\usepackage[section]{placeins} 
\input{preamble}

\begin{document}

\title[mode=title]{On the anti-aliasing properties of entropy filtering for discontinuous spectral element approximations of under-resolved turbulent flows}
\shorttitle{Anti-aliasing properties of entropy filtering}

\shortauthors{T. Dzanic \textit{et al.}}
\author[1]{T. Dzanic}[orcid=0000-0003-3791-1134]
\cormark[1]
\cortext[cor1]{Corresponding author}
\ead{tdzanic@tamu.edu}
\author[2]{W. Trojak}[orcid=0000-0002-4407-8956]
\author[1]{F. D. Witherden}[orcid=0000-0003-2343-412X]
\address[1]{Department of Ocean Engineering, Texas A\&M University, College Station, TX 77843, United States of America}
\address[2]{IBM Research, The Hartree Centre, Daresbury, WA4 4AD, United Kingdom}


\begin{abstract}
For large Reynolds number flows, it is typically necessary to perform simulations that are under-resolved with respect to the underlying flow physics. For nodal discontinuous spectral element approximations of these under-resolved flows, the collocation projection of the nonlinear flux can introduce aliasing errors which can result in numerical instabilities. In Dzanic and Witherden (\textit{J. Comput. Phys.}, 468, 2022), an entropy-based adaptive filtering approach was introduced as a robust, parameter-free shock-capturing method for discontinuous spectral element methods. This work explores the ability of entropy filtering for mitigating aliasing-driven instabilities in the simulation of under-resolved turbulent flows through high-order implicit large eddy simulations of a NACA0021 airfoil in deep stall at a Reynolds number of 270,000. It was observed that entropy filtering can adequately mitigate aliasing-driven instabilities without degrading the accuracy of the underlying high-order scheme on par with standard anti-aliasing methods such as over-integration, albeit with marginally worse performance at higher approximation orders.
\end{abstract}




\maketitle

\input{introduction}

\input{methodology}
\input{results}
\input{conclusion}

\bibliographystyle{unsrtnat}
\bibliography{reference}

\end{document}

%% file: preamble.tex
\usepackage{amssymb}
\usepackage{amsthm}
\usepackage{amsmath}
\usepackage{upgreek}
\usepackage{pdflscape}
\usepackage{listings}
\usepackage{multirow}
\usepackage[percent]{overpic}
\usepackage{multicol}
\usepackage{color}
\usepackage{stmaryrd}
\usepackage{xfrac}
\usepackage[capitalise]{cleveref}
\usepackage{booktabs}
\usepackage[section]{placeins} 
\usepackage[section]{algorithm}
\usepackage{calc}
\usepackage[titletoc]{appendix}
\usepackage{siunitx}
\usepackage[sort&compress,square,numbers]{natbib}

\usepackage{graphicx}
\usepackage{graphics}
\usepackage{wrapfig}
\usepackage{float}
\usepackage{subfig}
\usepackage[percent]{overpic}
\usepackage{varwidth}
\usepackage{tikz}
\usepackage{pgfplots}
\usepackage{adjustbox}
\usetikzlibrary{arrows,matrix,positioning,fit}
\usetikzlibrary{shapes,positioning}
\usetikzlibrary{backgrounds}
\usepackage{tikz-layers}
\usepgfplotslibrary{fillbetween}
\usetikzlibrary{intersections}
\pgfplotsset{compat=1.14}
\usepgfplotslibrary{colorbrewer}
\usepgfplotslibrary{patchplots}
\usepgfplotslibrary[colorbrewer]
\usetikzlibrary{pgfplots.colorbrewer}
\usetikzlibrary[pgfplots.colorbrewer]
\usepgfplotslibrary{units}
\usetikzlibrary{spy}
\usepackage{pgfplotstable}
\usepackage{arrayjobx}
\usepackage{pifont}
\graphicspath{ {./figs/} }
\usetikzlibrary{external}
\usetikzlibrary{patterns}

\newlength\myheight
\newlength\mydepth
\settototalheight\myheight{Xygp}
\newcommand*\inlinegraphics[1]{%
  \settototalheight\myheight{Xygp}%
  \settodepth\mydepth{Xygp}%
  \raisebox{-\mydepth}{\includegraphics[height=\myheight]{#1}}%
}
\newcommand\orcid[1]{\href{https://orcid.org/#1}{\inlinegraphics{orcid_16x16.png}}}

\makeatletter
\def\BState{\State\hskip-\ALG@thistlm}
\makeatother

\newdefinition{definition}{Definition}[section]

%% file: introduction.tex
\section{Introduction}
\label{sec:intro}
The accurate and efficient prediction of complex turbulent flows has been a significant driving force in the development of computational fluid dynamics methods over the decades. The primary difficulty in the simulation of these flows is the large variation of scales encountered in high Reynolds number regimes, such that the numerical resolution requirements of the underlying physical phenomena make direct simulation computationally intractable for practical applications. In many cases, it is necessary to perform simulations that are under-resolved with respect to the physical scales of the flow field in question. While this lack of resolution introduces error in the approximation of the flow, the goal of these under-resolved simulations is to still accurately predict the predominant large-scale flow physics at a lower computational cost. 

For these scale-resolving simulations of turbulent flows, discontinuous spectral element methods (DSEM) have shown promise as an efficient and accurate numerical approach~\citep{Hesthaven2008a, Huynh2007, Kopriva1996, Liu2004}. By approximating the solution by a set of piecewise-continuous high-order polynomials within each element, these schemes combine the geometric flexibility of finite volume methods with the arbitrarily high-order accuracy and efficiency of finite difference methods. As such, they have allowed for the computation of complex fluid flows that would otherwise be intractable with standard low-order numerical methods. However, for nonlinear equations such as the ones governing fluid flow, the projection of the nonlinear flux onto the polynomial space spanned by the solution can introduce numerical instabilities~\citep{Karniadakis2005}. These instabilities, so-called aliasing errors, stem from the energy of high-frequency under-resolved modes aliasing to lower frequency modes~\citep{Spiegel2015}. This effect becomes particularly problematic for simulations that are significantly under-resolved, where aliasing errors can result in unpredictable behavior or the failure of the numerical approach altogether~\citep{Gassner2012}. 

Various approaches have been proposed as numerical stabilization and anti-aliasing methods for under-resolved flows, including modal filtering~\citep{Hesthaven2008, Gassner2012, Gottlieb1997}, over-integration~\citep{Spiegel2015, Gassner2012}, and spectral vanishing viscosity~\citep{Tadmor1990,Karamanos2000}. Of these methods, the modal filtering approach is typically the simplest and most common technique, where the high-frequency content in the solution is explicitly dissipated to attempt to mitigate aliasing-driven instabilities. However, this simplicity comes with a lack of robustness, requiring problem-specific parameters that can require extensive tuning and may drastically affect the stability and accuracy of the simulation. On the other end of the spectrum, over-integration, sometimes referred to as polynomial de-aliasing or Galerkin projection, is a robust and accurate technique for mitigating aliasing errors. In this approach, the nonlinear flux is approximated using a more optimal projection onto the polynomial space spanned by the solution, typically performed by utilizing a sufficiently strong quadrature rule. While this method can offer superior accuracy in comparison to other approaches, it comes with an increase in computational cost which may be significant for severely under-resolved flows. As such, there is potential to expand the ability to robustly and efficiently simulate turbulent flows using DSEM if alternate approaches can be used to mitigate aliasing errors. 

In \citet{Dzanic2022}, an entropy-based adaptive filtering approach was introduced for the purpose of mitigating numerical instabilities stemming from high-order DSEM approximations of discontinuous flow features (i.e., shocks). It was observed by the authors that this shock-capturing approach, referred to as entropy filtering, also allowed for the simulation of high Reynolds number flows on under-resolved meshes that would typically be unstable due to aliasing errors. We posit that this is a result of aliasing errors manifesting as violations of a local minimum entropy principle, which is supported by the observations of \citet{Honein2004}, such that entropy filtering may adequately perform as an anti-aliasing technique separately from its purpose as a shock-capturing approach. In this work, we explore the capability and accuracy of entropy filtering for anti-aliasing in DSEM approximations of under-resolved turbulent flows. In particular, we consider implicit large eddy simulations of a NACA0021 in deep stall from the DESider project \citep{Haase2009} as presented by \citet{Park2017}, a case notorious for aliasing driven instabilities in high-order methods that requires a substantial amount of numerical stabilization for the given setup, and present a comparison to standard anti-aliasing approaches. 

The remainder of this paper is organized as follows. We present the methodology in \cref{sec:methodology}, including an overview of the numerical approach, various anti-aliasing techniques, and the problem setup. The results of the numerical experiments are then shown in \cref{sec:results}. Conclusions are finally drawn in \cref{sec:conclusion}.

%% file: methodology.tex
\section{Methodology}\label{sec:methodology}
\subsection{Governing equations}
We consider the three-dimensional compressible Navier--Stokes equations, given as
\begin{align}
    \frac{\partial \rho}{\partial t} + \frac{\partial}{\partial x_i}(\rho v_i) &= 0,\\
    \frac{\partial}{\partial t}\rho v_i + \frac{\partial}{\partial x_j}(\rho v_i v_j) &= -\frac{\partial P}{\partial x_i} + \frac{\partial}{\partial x_j}\tau_{ij},\\
    \frac{\partial}{\partial t}E + \frac{\partial}{\partial x_j}( E v_j) &= \frac{\partial}{\partial x_j}\left [v_i \tau_{ij} - q_j - P v_j\right],
\end{align}
where $\rho$ is the density, $\rho v_i$ are the momentum components, and $E$ is the total energy. Furthermore, we define the stress tensor, strain tensor, and heat flux as
\begin{equation}
    \tau_{ij} = 2 \mu \bigg [ S_{ij} - \frac{1}{3} \frac{\partial v_k}{\partial x_k} \bigg ],
\end{equation}
\begin{equation}
    S_{ij} = \frac{1}{2} \bigg ( \frac{\partial v_i}{\partial x_j} + \frac{\partial v_j}{\partial x_i} \bigg ),
\end{equation}
and
\begin{equation}
    q_j = \frac{\mu}{Pr} \frac{\partial h}{\partial x_j},
\end{equation}
respectively, where $\mu$ is the dynamic viscosity and $Pr = 0.71$ is the molecular Prandtl number, and $h$ is the specific enthalpy. Using the ideal gas equation of state, the pressure is then defined as
\begin{equation}
    P = (\gamma - 1)(E - \frac{1}{2}\rho v_i v_i),
\end{equation}
where $\gamma = 1.4$ is the specific heat ratio. 

The Navier--Stokes equations can also be conveniently represented in the form of a system of conservation laws as
\begin{equation}
    \partial_t \mathbf{u} + \boldsymbol{\nabla}{\cdot}\left(\mathbf{F}_I(\mathbf{u}) + \mathbf{F}_V(\mathbf{u}) \right) = 0,
\end{equation}
where $\mathbf{u} = [\rho, \boldsymbol{\rho}\mathbf{v}, E]^T = [\rho, \rho u, \rho v, \rho w, E]^T$ is the solution and
\begin{equation}\label{eq:navierstokes} \mathbf{F}_I = \begin{bmatrix}
            \boldsymbol{\rho}\mathbf{v} \\
            \boldsymbol{\rho}\mathbf{v} \otimes\mathbf{v}  + P\mathbf{I}\\
        (E+P)\mathbf{v}
    \end{bmatrix} \quad \mathrm{and} \quad \mathbf{F}_V = \begin{bmatrix}
            0\\
            -\mu \left(\nabla \mathbf{v} + \nabla \mathbf{v}^T \right)  + \frac{2}{3}\mu \nabla\cdot \mathbf{v} \\
        -\mu \left(\nabla \mathbf{v} + \nabla \mathbf{v}^T \right)\mathbf{v} -  \frac{\mu}{Pr} \nabla h 
    \end{bmatrix}
\end{equation}
are the inviscid and viscous fluxes, respectively.

\subsection{Numerical discretization}
The governing equations were discretized using the flux reconstruction scheme of \citet{Huynh2007}, a generalization of the nodal discontinuous Galerkin method \citep{Hesthaven2008a}. In this approach, the mesh $\mathcal T$ is a discretization of the domain $\Omega$ with $N$ disjoint elements, such that $\mathcal T = \bigcup_{N}\mathcal T_k$ and $\mathcal T_i\cap\mathcal T_j=\emptyset$ for $i\neq j$. Within each element $\mathcal T_k$, a discrete solution $\mathbf{u}_h(\mathbf{x})$ is obtained through a nodal interpolating approximation as
\begin{equation}
     \mathbf{u}_h (\mathbf{x}) = \sum_{i = 1}^{N_s} \mathbf{u} (\mathbf{x}^s_i) {\phi}_i (\mathbf{x}),
\end{equation}
where $\mathbf{x}^s_i \in \mathcal T_k \ \forall \ i \in \{1,..., N_s\}$ is a set of $N_s$ solution nodes and ${\phi}_i (\mathbf{x})$ is a set of polynomial basis functions with the property ${\phi}_i (\mathbf{x}^s_j) = \delta_{ij}$. We use the notation $\mathbb P_p$ to represent the order of the approximation for some order $p$, defined as the maximal order of $\mathbf{u}_h(\mathbf{x})$. 

A discontinuous approximation of the inviscid flux is first formed via a collocation projection of the inviscid flux onto the solution space, i.e., 
\begin{equation}
    \mathbf{f}^D (\mathbf{x}) = \sum_{i = 1}^{N_s} \mathbf{F}_I \left ( \mathbf{u} (\mathbf{x}^s_i) \right )  {\phi}_i (\mathbf{x}),
\end{equation}
A corrected flux is then formed by amending the discontinuous flux with additional correction terms which enforce $C^0$ continuity in the normal direction of $\partial\mathcal T_k$ as
\begin{equation}
    \mathbf{f}^C(\mathbf{x}) = \mathbf{f}^D(\mathbf{x}) + \sum_{i = 1}^{N_f} \left[\overline{\mathbf{F}}_i - \mathbf{f}^D(\mathbf{x}^f_i)\cdot \mathbf{n}_i \right] \mathbf{h}_i (\mathbf{x}),
\end{equation}
where $\mathbf{x}^f_i \in \partial \Omega_k \ \forall \ i \in \{1,..., N_f\}$ is a set of $N_f$ interface flux nodes, $\mathbf{n}_i$ is their associated outward-facing normal vector, $\overline{\mathbf{F}}_i$ is the common inviscid interface flux to be defined in \cref{ssec:setup}, and $\mathbf{h}_i$ is the correction function associated with the given flux node. These correction functions have the properties that
\begin{equation}
     \sum_{i = 1}^{N_f} \mathbf{h}_i (\mathbf{x}) \in RT_p \quad \mathrm{and} \quad \mathbf{n}_i \cdot\mathbf{h}_j (\mathbf{x}^f_i) = \delta_{ij},
\end{equation}
where $RT_p$ is the Raviart--Thomas space of order $p$ \citep{Raviart1977}. The correction functions are chosen to recover the nodal discontinuous Galerkin approach \citep{Huynh2007, Hesthaven2008, Trojak2021} in this work. 

For the viscous component, it is necessary first to form an appropriate approximation of the gradient of the solution, represented as
\begin{equation}
    \mathbf{w} \approx \nabla \mathbf{u}.
\end{equation}
Similarly to the calculation of the inviscid flux, a $C^0$ continuous approximation of the solution gradient is formed by amending the discontinuous approximation of the solution gradient as
\begin{equation}
    \mathbf{w}^C(\mathbf{x}) = \mathbf{w}^D(\mathbf{x}) + \sum_{i = 1}^{N_f} \left[\overline{\mathbf{u}}_i - \mathbf{u}(\mathbf{x}^f_i) \right] \nabla \mathbf{h}_i (\mathbf{x}),
\end{equation}
where $\overline{\mathbf{u}}_i$ is the common interface solution and 
\begin{equation}
    \mathbf{w}^D (\mathbf{x}) = \sum_{i = 1}^{N_s}  \mathbf{u} (\mathbf{x}^s_i) \nabla {\phi}_i (\mathbf{x}).
\end{equation}
The discontinuous approximation of the viscous flux can then be computed as
\begin{equation}
    \mathbf{g}^D (\mathbf{x}) = \sum_{i = 1}^{N_s} \mathbf{F}_V \left ( \mathbf{u} (\mathbf{x}^s_i), \mathbf{w}^C (\mathbf{x}^s_i)\right )  {\phi}_i (\mathbf{x}),
\end{equation}
after which adding the correction terms yields
\begin{equation}
    \mathbf{g}^C(\mathbf{x}) = \mathbf{g}^D(\mathbf{x}) + \sum_{i = 1}^{N_f} \left[\overline{\mathbf{G}}_i - \mathbf{g}^D(\mathbf{x}^f_i)\cdot \mathbf{n}_i \right] \mathbf{h}_i (\mathbf{x}).
\end{equation}
Similarly to the common interface solution, a common viscous interface flux $\overline{\mathbf{G}}_i$, which will be defined in \cref{ssec:setup}, is used at the interfaces. With this discretization, the temporal derivative of the solution can then be approximated as
\begin{equation}
    \partial_t \mathbf{u} = -\boldsymbol{\nabla}{\cdot} \mathbf{f}^C(\mathbf{x}) - \boldsymbol{\nabla}{\cdot} \mathbf{g}^C(\mathbf{x}),
\end{equation}
after which the solution can be advanced using a suitable temporal integration scheme. 

\subsection{Anti-aliasing techniques}
To present an evaluation of the entropy filtering approach as an anti-aliasing technique, we present a comparison of the approach to two standard anti-aliasing methods: over-integration and modal filtering. For brevity, we occasionally use the shorthand notation EF, OI, and MF to denote entropy filtering, over-integration, and modal filtering, respectively. 

For the over-integration approach, the goal is to find the optimal approximation of the analytic flux within the span of the solution space. In the $L^2$ norm, the best possible approximation is the $L^2$ projection, or Galerkin projection, of the flux onto the span of the solution space, i.e.,
\begin{equation}
    \mathbf{f}_{OI}(\mathbf{x}) = \underset{\mathbf{f}}{\arg \min} \ \int_{\mathcal T_k} \| \mathbf{F}\left (\mathbf{u}_h(\mathbf{x}) \right) - \mathbf{f}(\mathbf{x}) \|_2 \ \mathrm{d}\mathbf{x}.
\end{equation}
The projected flux polynomial can be represented in modal form as
\begin{equation}
    \mathbf{f}_{OI}(\mathbf{x}) = \sum_{i=0}^{N_s} \hat{\mathbf{f}}_i {\psi}_i (\mathbf{x}),
\end{equation}
where ${\psi}_i (\mathbf{x})$ are a set of  modal basis functions that are orthogonal with respect to the unit measure and $\hat{\mathbf{f}}_i$ are their respective coefficients computed by the integral 
\begin{equation}
    \hat{\mathbf{f}}_i = \int_{\mathcal T_k}\mathbf{F}\left (\mathbf{u}_h(\mathbf{x}) \right){\psi}_i (\mathbf{x})\ \mathrm{d}\mathbf{x}.
\end{equation}
For non-polynomial flux functions, this integration usually cannot be performed exactly, but it may be approximated utilizing a suitably-strong quadrature rule, i.e.,
\begin{equation}
    \int_{\mathcal T_k}\mathbf{F}\left (\mathbf{u}_h(\mathbf{x}) \right){\psi}_i (\mathbf{x})\ \mathrm{d}\mathbf{x} \approx \sum_{i=j}^{N_q}w_j\mathbf{F}\left (\mathbf{u}_h(\mathbf{x}_j^q) \right){\psi}_i (\mathbf{x}_j^q),
\end{equation}
where $\mathbf{x}_j^q$ are a set of $N_q$ quadrature nodes (with $N_q \geq N_s$) and $w_j$ are their associated quadrature weights. Typically, the more under-resolved the flow is, the stronger the quadrature rule has to be to mitigate aliasing errors, such that for highly under-resolved flows, the cost of computing this projection may become substantial. In this work, we utilize Gaussian quadrature rules to compute the projection with the notation OI-$\mathbb Q_q$ referring to over-integration with a $q$-th degree quadrature rule. For a more in-depth description of this approach, the reader is referred to \citet{Park2017}, Section II.D.

For the modal filtering approach, the goal is to mitigate aliasing driven instabilities by explicitly filtering the high-frequency content from the solution as these errors tend to manifest at higher frequencies \citep{Cox2021}. In this approach, the solution is first transformed to its modal form $\widetilde{\mathbf{u}}_h(\mathbf{x})$, defined as
\begin{equation}
    \widehat{\mathbf{u}}_h(\mathbf{x}) = \sum_{i = 1}^{N_s} \widehat{\mathbf{u}}_i {\psi}_i (\mathbf{x}) = \mathbf{u}_h(\mathbf{x}),
\end{equation}
where $\widehat{\mathbf{u}}_i$ are the associated modes corresponding to the modal basis functions ${\psi}_i (\mathbf{x})$. In this modal form, a filtered solution $\widetilde{\mathbf{u}}(\mathbf{x})$ can then be obtained by applying a filtering operation $H$to the individual modes as
\begin{equation}
    \widetilde{\mathbf{u}}_h(\mathbf{x}) = \sum_{i = 1}^{N_s} H_i \left( \widehat{\mathbf{u}}_i \right) {\psi}_i (\mathbf{x})
\end{equation}
A standard choice of filter kernel is the exponential filter \citep{Hesthaven2008}, given as
\begin{equation}
    H_i(\widehat{\mathbf{u}}_i) =
    \begin{cases}
    \widehat{\mathbf{u}}_i, \quad \quad \quad \quad \quad \quad \quad \quad \quad \quad \, \mathrm{if} \ \eta_i \leq \eta_c, \\
    \widehat{\mathbf{u}}_i\exp \left [- \kappa \left(\frac{\eta_i - \eta_c}{\eta_m - \eta_c}\right)^s \right], \quad \quad \mathrm{else},
    \end{cases}
\end{equation}
where $\kappa \approx \log \epsilon$ for some value of machine precision $\epsilon$, $s$ is some even integer representing the filter order, $\eta_i$ is the maximal order of the $i$-th basis function, $\eta_c$ is the cutoff order, and $\eta_m = p+1$ is the maximal order. The optimal parameters for the modal filtering approach are typically not known \textit{a priori} and must be tuned on a per-case basis to yield an accurate and robust stabilization method.

For the entropy filtering approach, a filtering kernel is applied similarly to a modal filter, but its parameters are computed adaptively based on the solution's ability to preserve certain invariants of the system such as the positivity of density and pressure and a local discrete minimum entropy principle \citep{Dzanic2022}. This adaptive method is implemented in the context of a second-order exponential filter applied at each stage of the temporal scheme, given as
\begin{equation}
    H_i(\widehat{\mathbf{u}}_i) =
    \hat{\mathbf{u}}_i\exp \left (-\zeta \eta_i^2 \right),
\end{equation}
where $\zeta$ is the filter strength. This filter strength is computed via an element-wise nonlinear optimization problem that finds the minimum necessary filter strength such that the discrete nodal solution values have positive density ($\rho > 0$), positive pressure ($P > 0$), and an entropy greater than some minimum threshold ($\sigma > \sigma_{\min}$), i.e.,
\begin{equation}
    \zeta = \underset{\zeta\ \geq\ 0}{\mathrm{arg\ min}} \ \ \mathrm{s.t.} \ \  \left [\widetilde{\rho}_i \geq 0, \ \widetilde{P}_i \geq 0, \ \widetilde{\sigma}_i \geq \sigma_{\min} \ \ \forall \ i \in \{1,...,N_s\}\right ].
\end{equation}
The entropy functional is taken as the specific physical entropy $\sigma = P\rho^{-\gamma}$, and the minimum entropy is computed as the minimum discrete entropy within the element and its Voronoi neighbors. For an in-depth overview of this approach and description of a computationally efficient implementation thereof, the reader is referred to \citet{Dzanic2022}, Section 3 and \citet{Dzanic2023}, Section 4.1. 

\subsection{Problem setup and computational framework}\label{ssec:setup}

The problem setup consists of a NACA0021 airfoil operating at a Reynolds number of $Re = 270,000$, Mach number of $M = 0.1$, and an angle of attack $\alpha = 60^{\circ}$. At these operating conditions, the airfoil is in deep stall, with strongly separated flow on the suction side of the wing which yields complex unsteady flow physics. If these physics are not very well-resolved, aliasing errors can quickly cause the simulation to diverge. For a suitable comparison of the anti-aliasing capabilities of the three techniques, we utilize an identical computational setup as the work of \citet{Park2017}, for which the mesh resolution and approximation orders were chosen such as to cause significant numerical instabilities without sufficient anti-aliasing. 

\begin{figure}[tbhp]
    \centering
    \subfloat[$\mathbb P_3$]{\adjustbox{width=0.48\linewidth, valign=b}{\includegraphics[]{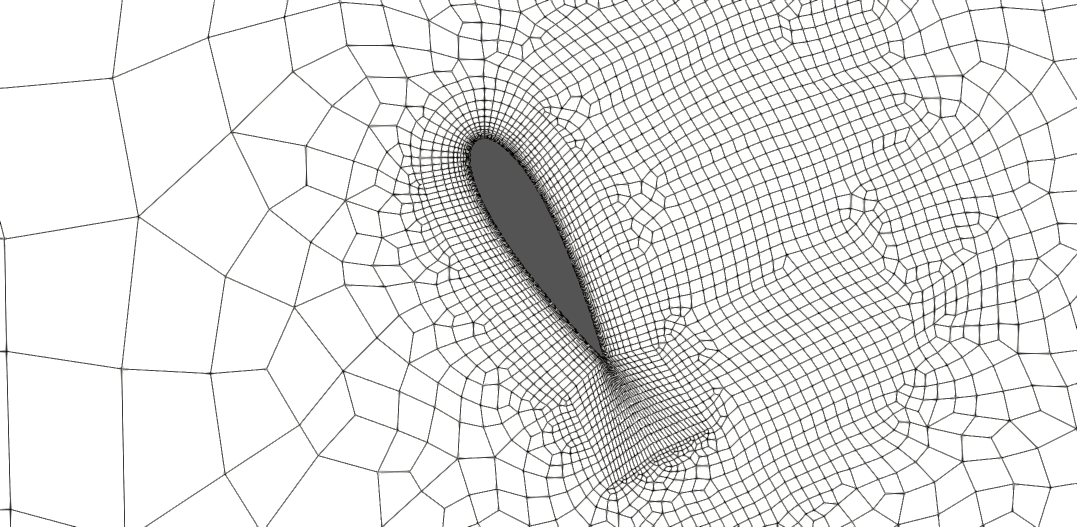}}}
    \hspace{10pt}
    \subfloat[$\mathbb P_4$]{\adjustbox{width=0.48\linewidth, valign=b}{\includegraphics[]{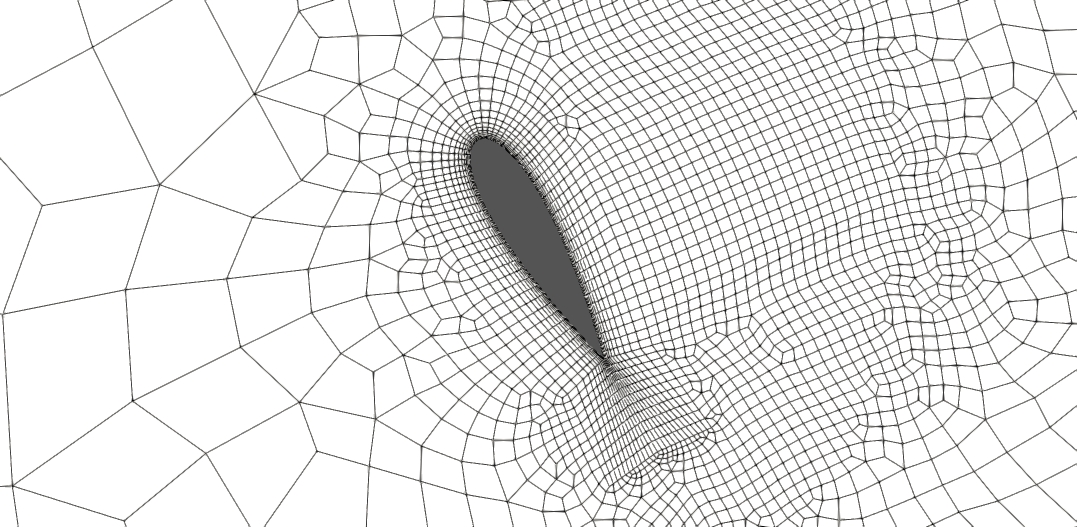}}}
    \caption{\label{fig:mesh} Cross-section of the mesh used for the $\mathbb P_3$ (left) and $\mathbb P_4$ (right) simulations.}
\end{figure}

The problem was solved using both a $\mathbb P_3$ and $\mathbb P_4$ approximation, corresponding to nominally fourth-order and fifth-order accurate spatial discretizations, respectively. The unstructured hexahedral meshes of \citet{Park2017} were used, shown in \cref{fig:mesh}, consisting of a finer $\mathbb P_3$ mesh (323,360 hexahedral elements) and a coarser $\mathbb P_4$ mesh (206,528 hexahedral elements). The $\mathbb P_4$ mesh was appropriately coarsened such that the degrees of freedom and the relative nodal spacing between the two approximation orders were roughly equal, and both meshes were generated such that the wall-normal spacing of the first solution point was $y^+ \approx 1$. An aspect ratio of 4 was chosen as this was found to be sufficiently large enough to mitigate the sensitivity of the flow to the spanwise extent \citep{Park2017}.

Simulations were performed using PyFR \citep{Witherden2014}, a high-order flux reconstruction solver that can efficiently target massively-parallel CPU and GPU computing architectures. Computations were performed on 16 NVIDIA V100 GPUs. Common inviscid interface fluxes were computed using a Rusanov-type \citep{Rusanov1962} Riemann solver with the Davis wavespeed estimate \citep{Davis1988}, and the common viscous interface fluxes were computed with the BR2 approach of \citet{Bassi2005}. The solution nodes were placed at the Gauss--Legendre quadrature nodes for the over-integration and modal filtering approaches and the Gauss--Lobatto quadrature nodes for the entropy filtering approach, the latter of which is due to the requirement of collocated solution and interface flux points for the method. Temporal integration was performed using the classic fourth-order, four stage Runge--Kutta scheme with a fixed time step. \textcolor{red}{For consistency with the work of \citet{Park2017}}, the simulations were initialized with $\mathbb P_1$ at $Re = 27,000$ and run until a characteristic time of $t_c = c/U_{\infty} = 50$, where $c$ is the chord and $U_{\infty}$ is the freestream velocity. The Reynolds number and approximation order were then increased to the operating conditions and the anti-aliasing methods of choice were applied \textcolor{red}{as no continuation in approximation order was found to be necessary to mitigate any instabilities stemming from the initial transients in the solution}. Without any anti-aliasing approach, the simulation quickly diverged. The flow was allowed to develop until $t_c = 100$, and then averaging was performed over the range $t_c \in [100, 300]$. Increasing the averaging period was found to have a minimal effect on the results. 

For anti-aliasing, $\mathbb Q_9$ and $\mathbb Q_{11}$ over-integration was used for the $\mathbb P_3$ and $\mathbb P_4$ approximations, respectively, as this was deemed to be the necessary amount of anti-aliasing to stabilize the solution over the simulation time \citep{Park2017}. The modal filter parameters were chosen as $\kappa = 32$, $\eta_c = 0$, and $s = 8$ (for $\mathbb P_3$) and $s = 6$ (for $\mathbb P_4$), with the filter applied every $N = 20$ time steps. These values were found by systematically modulating each component to find a set of values that stabilized the solution over the simulation period. We remark here that it is unlikely that these are the optimal parameters for this approach, and there could be values of these parameters that yield a stable approach with more accurate results. However, this is the typical drawback of stabilization techniques that have free parameters---it can be extremely costly to attempt to optimize these parameters over many simulations, and, without a proper point of reference, it is ambiguous as to which parameters result in more accurate predictions.
Furthermore, due to the different stabilizing effects of the various anti-aliasing approaches, different maximum time steps were permissible. For the over-integration approach, the maximum allowable time step was $\Delta t = 2{\cdot}10^{-5}$, whereas the MF approach only allowed a time step of $\Delta t = 1{\cdot}10^{-5}$. However, due to its additional nonlinear stability properties, the EF approach allowed a larger time step of $\Delta t = 5{\cdot}10^{-5}$~\citep{Trojak2018}. These factors contribute to the overall cost of the simulations in addition to the per-time-step cost of the various anti-aliasing approaches. 

%% file: results.tex
\section{Results}\label{sec:results}
\subsection{Pressure distributions}
\begin{figure}[htbp!]
    \centering
    \subfloat[$\mathbb P_3$]{\adjustbox{width=0.48\linewidth, valign=b}{\input{figs/naca_p3_cp}}}
    \subfloat[$\mathbb P_4$]{\adjustbox{width=0.48\linewidth, valign=b}{\input{figs/naca_p4_cp}}}
    \caption{\label{fig:cp} Average surface pressure coefficient distribution computed using a $\mathbb P_3$ approximation (left) and $\mathbb P_4$ approximation (right) with over-integration (OI), modal filtering (MF), and entropy filtering (EF). Experimental results of \citet{Swalwell2005} shown for reference. }
\end{figure}
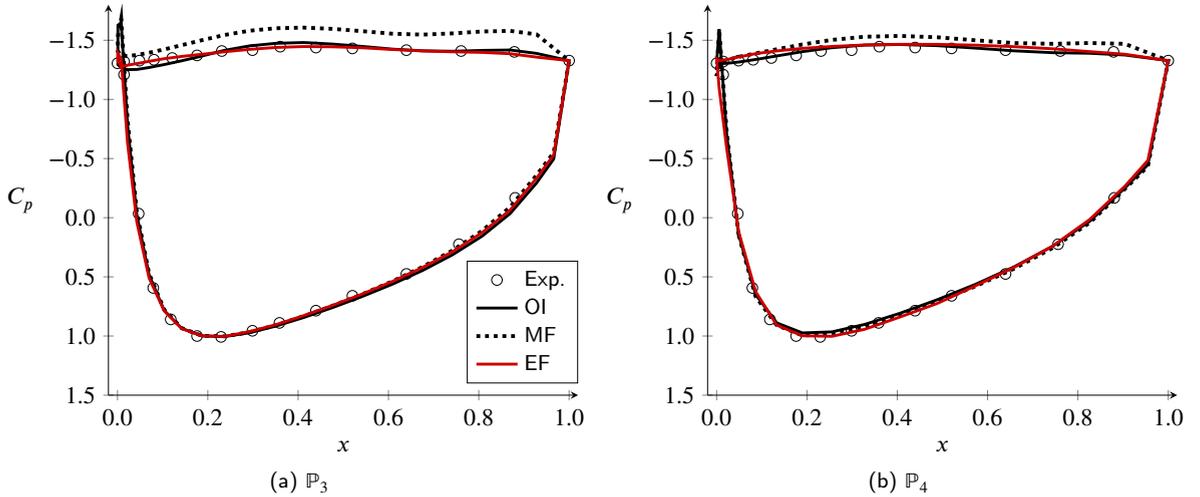

The results were first analyzed with respect to the average surface forces on the wing. The time- and span-averaged surface pressure coefficient distributions are shown in \cref{fig:cp} for the various anti-aliasing approaches computed with the $\mathbb P_3$ and $\mathbb P_4$ approximations. The experimental results of \citet{Swalwell2005} are additionally shown for reference. For both approximation orders, all anti-aliasing approaches showed good predictions of the pressure side surface pressure coefficient distribution, which is expected due to the attached flow in the region. On the suction side, where the highly unsteady separated region makes the flow physics much more complex, more drastic differences were observed. At $\mathbb P_3$, the over-integration approach showed good agreement with the experimental results, \textcolor{red}{with} only a marginal under-prediction of the (negative) pressure coefficient aft of the leading edge. In contrast, the modal filtering approach showed a large over-prediction of the (negative) pressure coefficient on the suction side, with a relatively consistent error across the length of the airfoil. In comparison, the entropy filtering approach showed predictions notably similar to the over-integration approach with even better agreement near the leading edge, although the differences between the two were quite marginal. 

When the approximation order was increased to $\mathbb P_4$, the variation in the results between the various approaches diminished. The over-integration approach showed excellent agreement with the reference data, marginally better than the over-integration results at $\mathbb P_3$. The modal filtering approach showed the largest difference in the results with increasing approximation order, with a notable improvement in the prediction of the suction side pressure coefficient distribution. However, a noticeable overprediction in the (negative) pressure coefficient was still observed, such that the approach showed relatively poor agreement with the experimental data. In contrast, the entropy filtering approach still showed good agreement with the over-integration results and the experimental results, with only a marginal overprediction of the (negative) pressure coefficient aft of the leading edge. Of the three different anti-aliasing techniques, the over-integration and modal filtering approaches showed improvements with increasing approximation order, with over-integration showing very marginal improvements and modal filtering showing noticeable improvements. In contrast, the $\mathbb P_3$ entropy filtering results were closer to experimental results and the $\mathbb P_4$ over-integration results, the latter of which can be considered as the reference numerical results for the given case setup, than the $\mathbb P_4$ results, although these differences were very minor. Interestingly enough, the $\mathbb P_3$ entropy filtering results were closer to the $\mathbb P_4$ over-integration than the $\mathbb P_3$ over-integration results. 

\subsection{Force coefficients}

The average lift and drag coefficients as computed by the various approaches were calculated, and the error in these coefficients was computed with respect to the experimental results of \citet{Swalwell2005} which report a lift coefficient of $C_L = 0.931$ and a drag coefficient of $C_D = 1.517$. The calculated force coefficients and errors are tabulated in \cref{tab:forces}. We remark here that the magnitude of the errors is not necessarily indicative of the accuracy of the anti-aliasing approach or the simulation, and errors of up to $20\%$ have been reported in various numerical experiments \citep{Garbaruk2003}. However, they are still presented as a point of comparison. At $\mathbb P_3$, the entropy filtering results showed the best agreement with the experimental data, with a $0.7\%$ error in the lift coefficient and a $3.0\%$ error in the drag coefficient. While the surface pressure coefficient distributions of over-integration approach and the entropy filtering approach were very similar, the over-integration approach showed a significantly larger error in the force coefficients than the entropy filtering approach, with an $8.4\%$ error in the lift coefficient and an $8.7\%$ error in the drag coefficient. This error was on par with the error from the modal filtering approach which was marginally higher. When the approximation order was increased to $\mathbb P_4$, the over-integration results showed the best agreement with the experimental data, with a $1.9\%$ error in the lift coefficient and a $3.0\%$ error in the drag coefficient. Consistent with the observations in the surface pressure coefficient distributions, the $\mathbb P_3$ entropy filtering results were most similar to the $\mathbb P_4$ over-integration results, and the modal filtering approach showed improvements in accuracy with increasing approximation order whereas the entropy filtering approach showed a degradation in accuracy.

\begin{table}[htbp!]
    \centering
    \begin{tabular}{lcccc}
    \toprule
    Method & $C_L$ &  $C_D$ & $C_L$ error &  $C_D$ error\\ 
    \midrule
    Experiment &  0.931 & 1.517 & - & -\\
    \midrule
    $\mathbb P_3$-OI &  1.009 & 1.650 & 8.4\% & 8.7\% \\
    $\mathbb P_3$-MF &  1.018 & 1.680 & 9.3\% & 10.8\% \\
    $\mathbb P_3$-EF &  0.937 & 1.562 & 0.7\% & 3.0\% \\
    \midrule
    $\mathbb P_4$-OI &  0.949 & 1.563 & 1.9\% & 3.0\% \\
    $\mathbb P_4$-MF &  0.977 & 1.625 & 4.9\% & 7.1\% \\
    $\mathbb P_4$-EF &  0.886 & 1.475 & 4.8\% & 2.7\% \\
    \bottomrule
    \end{tabular}
    \caption{\label{tab:forces}Average lift and drag coefficient computed using a $\mathbb P_3$ and $\mathbb P_4$ approximation with over-integration (OI), modal filtering (MF), and entropy filtering (EF). Error computed with respect to experimental results of \citet{Swalwell2005}. }
\end{table}

\subsection{Flow fields}
\begin{figure}[tbhp]
    \centering
    \subfloat[Over-integration]{\adjustbox{width=0.33\linewidth, valign=b}{\includegraphics[]{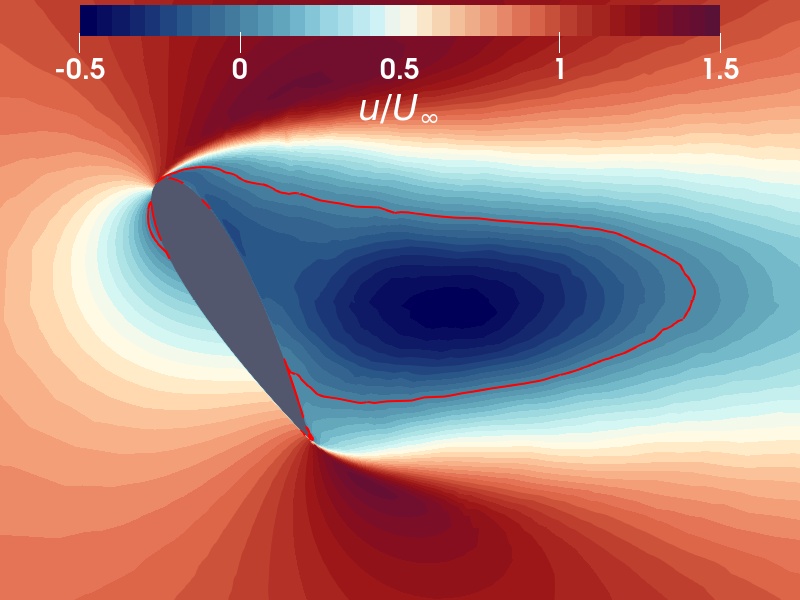}}}
    \hfill
    \subfloat[Modal filter]{\adjustbox{width=0.33\linewidth, valign=b}{\includegraphics[]{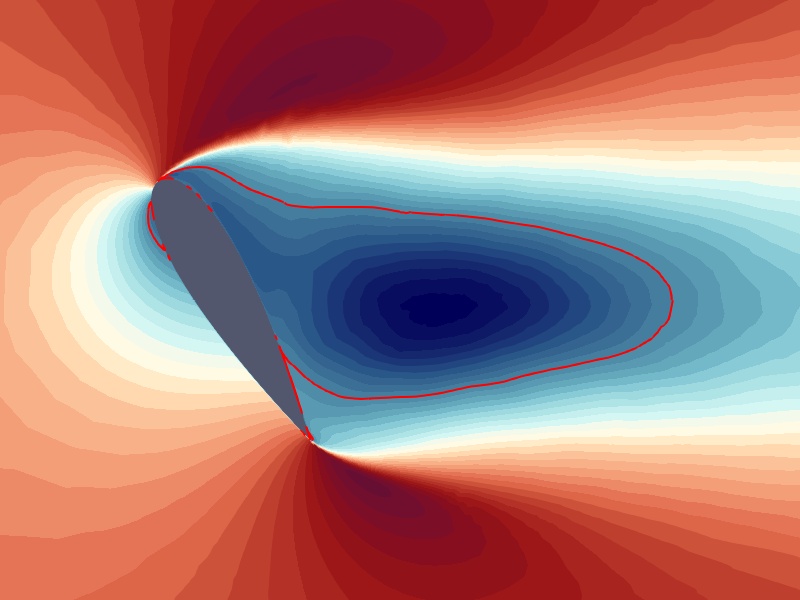}}}
    \hfill
    \subfloat[Entropy filter]{\adjustbox{width=0.33\linewidth, valign=b}{\includegraphics[]{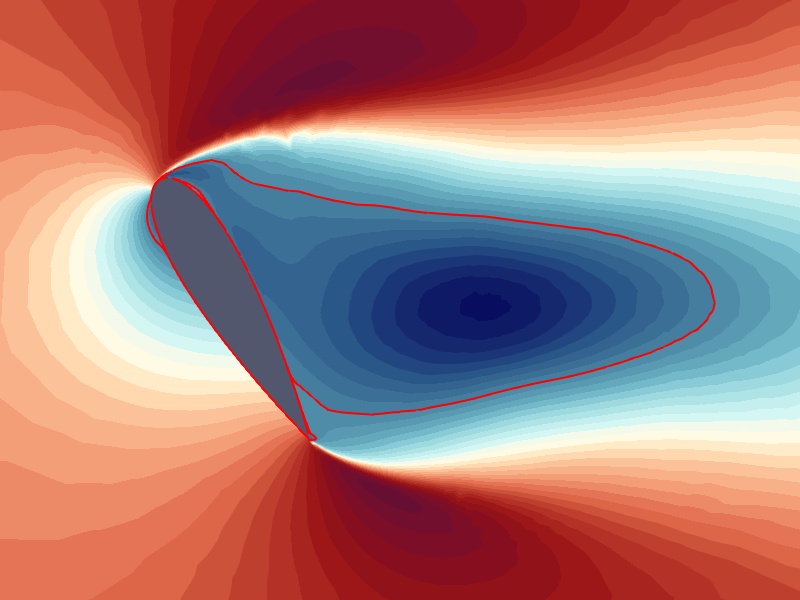}}}
    \caption{\label{fig:u_p3} Contours of streamwise velocity computed using a $\mathbb P_3$ approximation with over-integration (left), modal filtering (middle), and entropy filtering (right). Zero velocity represented by red isocontour. }
\end{figure}

The average flow characteristics were then analyzed with respect to the flow in the separation region on the suction side of the airfoil. The time- and span-averaged streamwise velocity contours are shown in \cref{fig:u_p3} as computed by the varying anti-aliasing approaches with a $\mathbb P_3$ approximation. It can be seen that all approaches predict a similar flow profile, with a large separation region downstream of the wing. However, the magnitude of the flow reversal and the shape of the separation bubble differed, the latter of which is represented by a red isocontour in the figures. The over-integration approach and the modal filtering approach gave very similar predictions of the shape of the separation region as well as the magnitude of flow reversal in the core of the separation bubble. In contrast, the entropy filtering approach showed some noticeable differences, with a lower degree of maximum flow reversal in the separation bubble as well as a marginally taller separation region. \textcolor{red}{Furthermore, the entropy filtering approach showed more oscillatory behavior along the shear layer in the wake, indicating that the approach is less dissipative but potentially more prone to predicting spurious instabilities in the flow, and it is unclear whether this behavior is expected to generalize to other problems.}

\begin{figure}[htbp!]
    \centering
    \subfloat[Over-integration]{\adjustbox{width=0.33\linewidth, valign=b}{\includegraphics[]{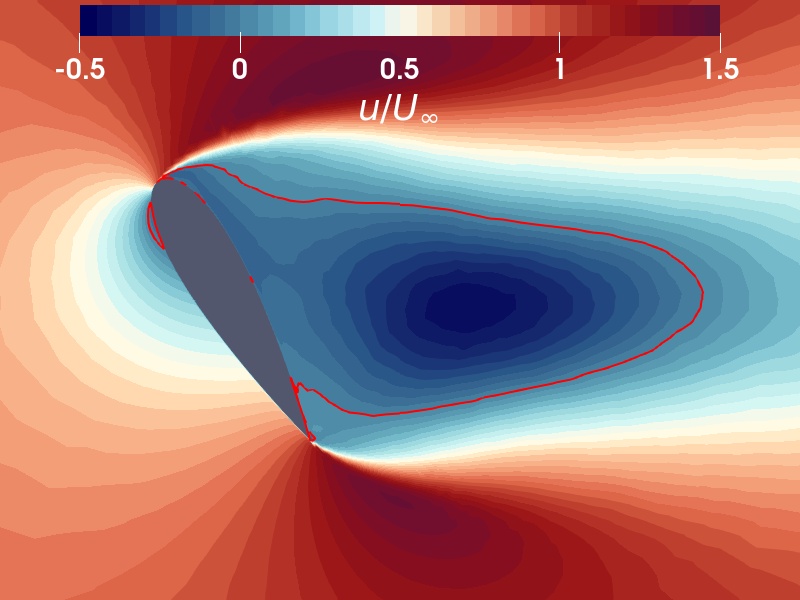}}}
    \hfill
    \subfloat[Modal filter]{\adjustbox{width=0.33\linewidth, valign=b}{\includegraphics[]{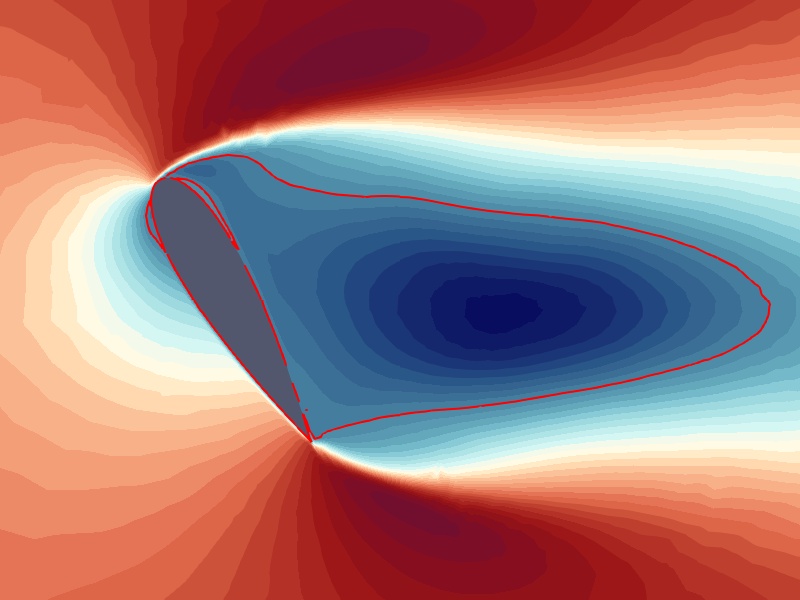}}}
    \hfill
    \subfloat[Entropy filter]{\adjustbox{width=0.33\linewidth, valign=b}{\includegraphics[]{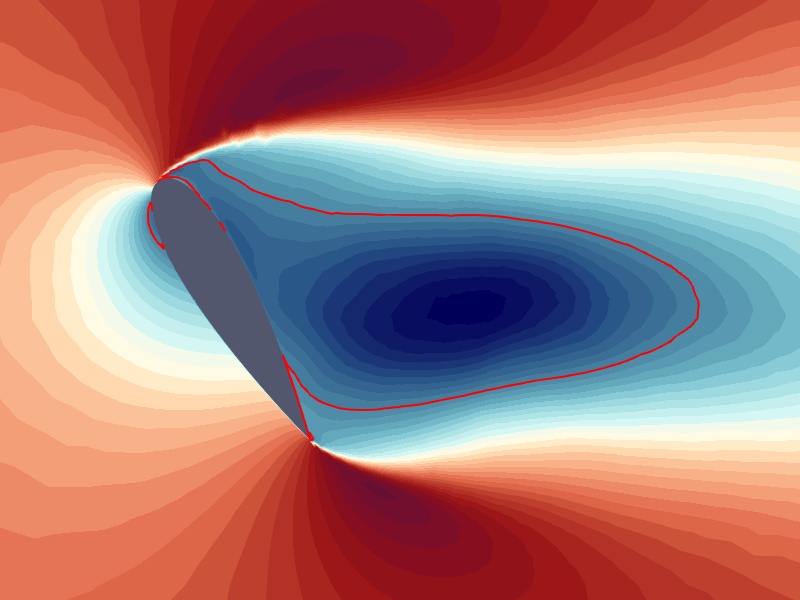}}}
    \caption{\label{fig:u_p4} Contours of streamwise velocity computed using a $\mathbb P_4$ approximation with over-integration (left), modal filtering (middle), and entropy filtering (right). Zero velocity represented by red isocontour. }
\end{figure}

Larger differences between the varying approaches were observed when the approximation order was increased to $\mathbb P_4$, shown by the contours in \cref{fig:u_p4}. With the over-integration approach, the magnitude of the maximum flow reversal decreased in comparison to the $\mathbb P_3$ results, while the shape of the separation region remained relatively similar. For the modal filtering approach, a large variation in the flow was seen when increasing the approximation order, with significant differences in the size and shape of the separation region as well as the location of the maximum magnitude of flow reversal. \textcolor{red}{This can likely be attributed to the sensitivity of the modal filtering approach to the choice of filtering parameters, the optimality of which may be dependent of the approximation order.} In comparison, the entropy filtering approach did not show as much of a change in the predicted flow field with the increase \textcolor{red}{in} the approximation order, with relatively similar predictions of the size and shape of the separation region. Similarly to the observations for the body forces and surface pressure coefficient distribution, the $\mathbb P_3$ entropy filtering results showed marked similarities to the $\mathbb P_4$ over-integration results, with nearly identical streamwise flow fields between the two approaches. \textcolor{red}{Furthermore, the oscillatory behavior along the shear layer in the wake in the entropy filtering approach diminished with increased approximation order, showing a similar shear layer profile as with the other two approaches. This is consistent with the observations in \citet{Dzanic2022} in that the entropy filtering approach as a shock capturing method was more dissipative at higher approximation orders. }

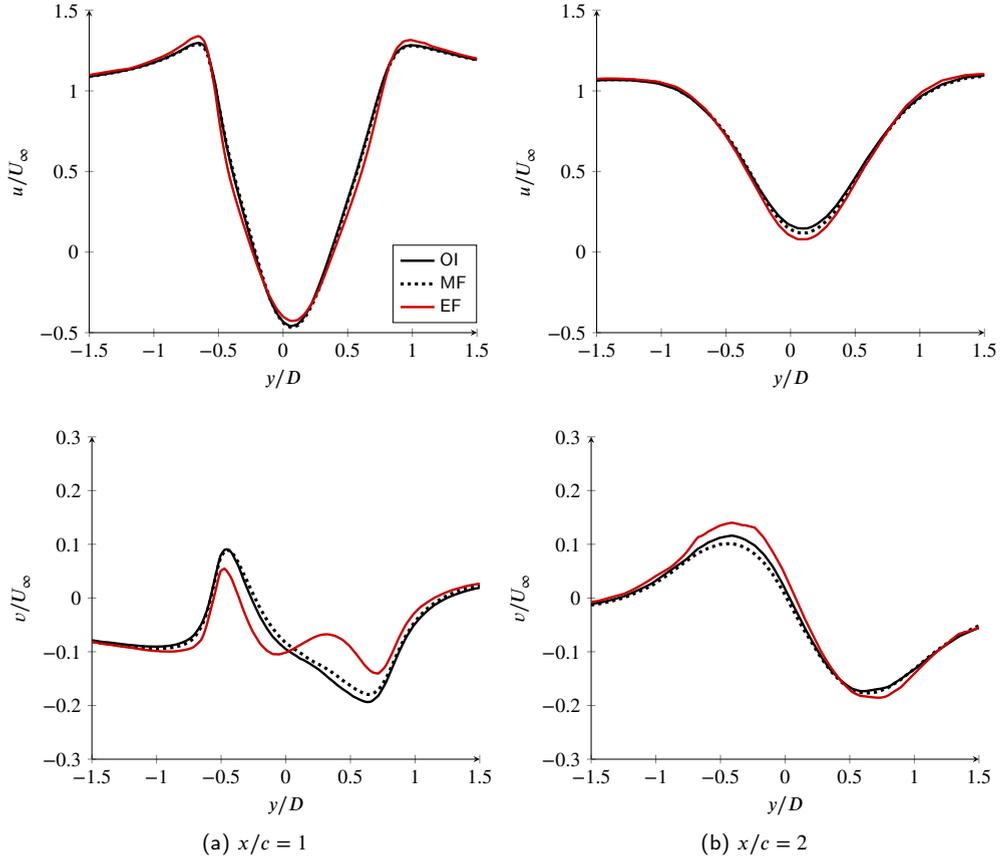
\begin{figure}[htbp!]
    \centering
    \adjustbox{width=0.4\linewidth, valign=b}{\input{figs/u_x1_p3}}
    \adjustbox{width=0.4\linewidth, valign=b}{\input{figs/u_x2_p3}}
    \newline
    \subfloat[$x/c = 1$]{\adjustbox{width=0.4\linewidth, valign=b}{\input{figs/v_x1_p3}}}
    \subfloat[$x/c = 2$]{\adjustbox{width=0.4\linewidth, valign=b}{\input{figs/v_x2_p3}}}
    \newline
    \caption{\label{fig:p3_profiles} Profiles of streamwise velocity (top row) and normal velocity (bottom row) at $x/c = 1$ (left) and $x/c = 2$ (right) computed using a $\mathbb P_3$ approximation with over-integration (OI), modal filtering (MF), and entropy filtering (EF).}
\end{figure}

The flow fields were further compared by analyzing the streamwise and normal velocity profiles in the wake. Two streamwise cross-section locations were chosen, one at $x/c = 1$, corresponding to the approximate location of the highest degree of flow reversal, and one at $x/c = 2$, corresponding to the approximate edge of the separation region. The time- and span-averaged streamwise and normal velocity profiles at these cross-section locations are shown in \cref{fig:p3_profiles} as computed by the varying anti-aliasing approaches with a $\mathbb P_3$ approximation. At both streamwise locations, the streamwise velocity profiles between the three approaches were similar, with only minor differences in the peak velocity defect. However, larger differences could be observed in the normal velocity profiles. Whereas the over-integration and modal filtering approaches showed very similar predictions, the entropy filtering approach predicted a distinct normal velocity profile in the separation region, showing a change in inflection near the centerline. Nevertheless, these differences diminished farther in the wake, with similar profiles between the three approaches at $x/c = 2$. 

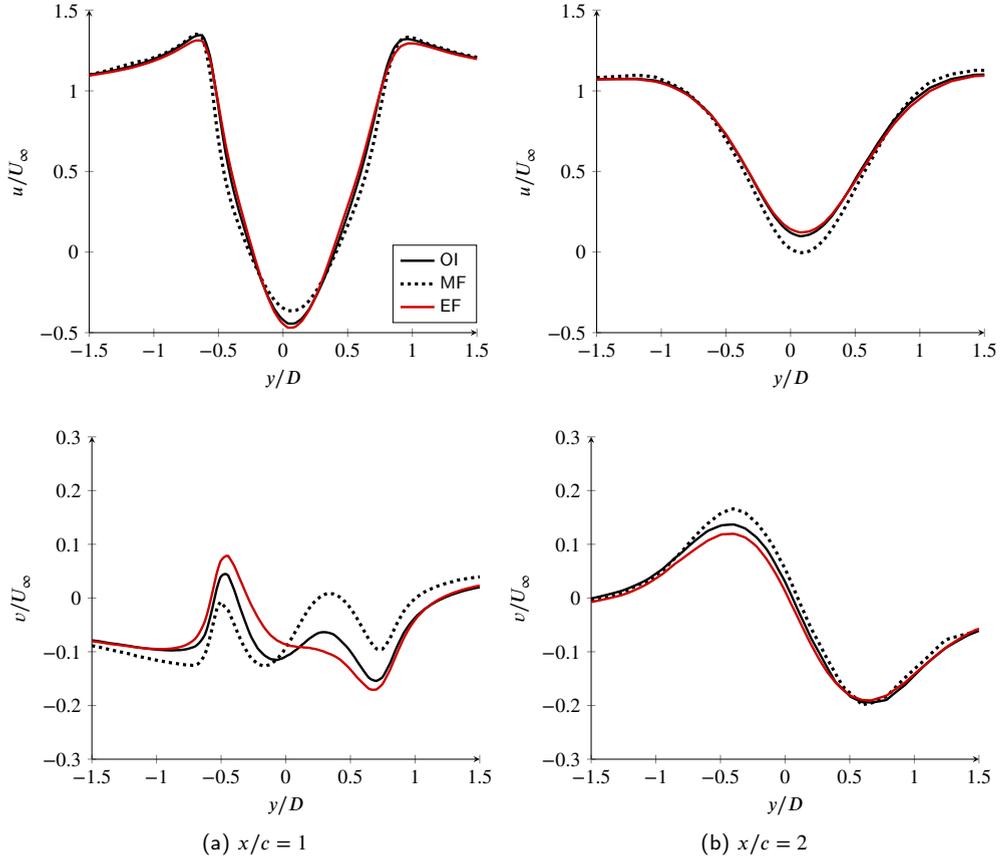
\begin{figure}[htbp!]
    \centering
    \adjustbox{width=0.4\linewidth, valign=b}{\input{figs/u_x1_p4}}
    \adjustbox{width=0.4\linewidth, valign=b}{\input{figs/u_x2_p4}}
    \newline
    \subfloat[$x/c = 1$]{\adjustbox{width=0.4\linewidth, valign=b}{\input{figs/v_x1_p4}}}
    \subfloat[$x/c = 2$]{\adjustbox{width=0.4\linewidth, valign=b}{\input{figs/v_x2_p4}}}
    \newline
    \caption{\label{fig:p4_profiles} Profiles of streamwise velocity (top row) and normal velocity (bottom row) at $x/c = 1$ (left) and $x/c = 2$ (right) computed using a $\mathbb P_4$ approximation with over-integration (OI), modal filtering (MF), and entropy filtering (EF).}
\end{figure}

With an increase in the approximation order to $\mathbb P_4$, more noticeable differences were observed between the three anti-aliasing approaches. The time- and span-averaged streamwise and normal velocity profiles as computed by the $\mathbb P_3$ approximation are shown in \cref{fig:p4_profiles}. In the streamwise velocity profiles, the over-integration and entropy filtering approaches showed very similar results, such that the profiles were essentially indistinguishable in both the separation region and farther in the wake. The modal filtering approach, however, showed distinct results in the streamwise profiles, with an underprediction of the velocity defect in the separation region and an overprediction of the velocity defect farther in the wake in comparison to the other two anti-aliasing approaches. \textcolor{red}{We remark here that these comparisons are strictly with respect to the differences between the various anti-aliasing approaches as there is a lack of a well-established ``true'' solution.} 
The differences between the three approaches were most pronounced in the normal velocity profiles in the separation region. The entropy filtering approach showed the most pronounced variation in the normal velocity profile, whereas the modal filtering approach showed the least. The over-integration approach showed a normal velocity profile very similar to the results of the $\mathbb P_3$ anti-aliasing approach, which is consistent with previous observations. These differences in the three approaches diminished farther in the wake, with the profiles at $x/c = 2$ showing very similar predictions.

\subsection{Force spectra}
The efficacy of the anti-aliasing approaches was further evaluated by analyzing the temporal statistics of the force coefficients. The power spectral density (PSD) of the lift coefficient was computed using Welch's averaged periodogram method with a sampling rate of $1/160 t_c$, window length of 4096, and a shift of 10. The PSD profiles of the various approaches were compared to the experimental results of \citet{Swalwell2005}, which show two distinct peaks in the PSD. We remark here that the experimental results were obtained using a sectional lift coefficient at a fixed spanwise location instead of the total lift coefficient as per the simulations. As such, the comparison is performed only with respect to the ability of predicting the frequencies of the dominant peaks in the PSD. 

\begin{figure}[htbp!]
    \centering
    \subfloat[$\mathbb P_3$]{\adjustbox{width=0.48\linewidth, valign=b}{\input{figs/naca_p3_psd}}}
    \subfloat[$\mathbb P_4$]{\adjustbox{width=0.48\linewidth, valign=b}{\input{figs/naca_p4_psd}}}
    \hfill
    \caption{\label{fig:psd} Lift coefficient power spectral density computed using a $\mathbb P_3$ approximation (left) and $\mathbb P_4$ approximation (right) with over-integration (OI), modal filtering (MF), and entropy filtering (EF). Experimental results of \citet{Swalwell2005} shown for reference. Primary and secondary PSD peaks in the experimental data shown as vertical dotted lines.}
\end{figure}
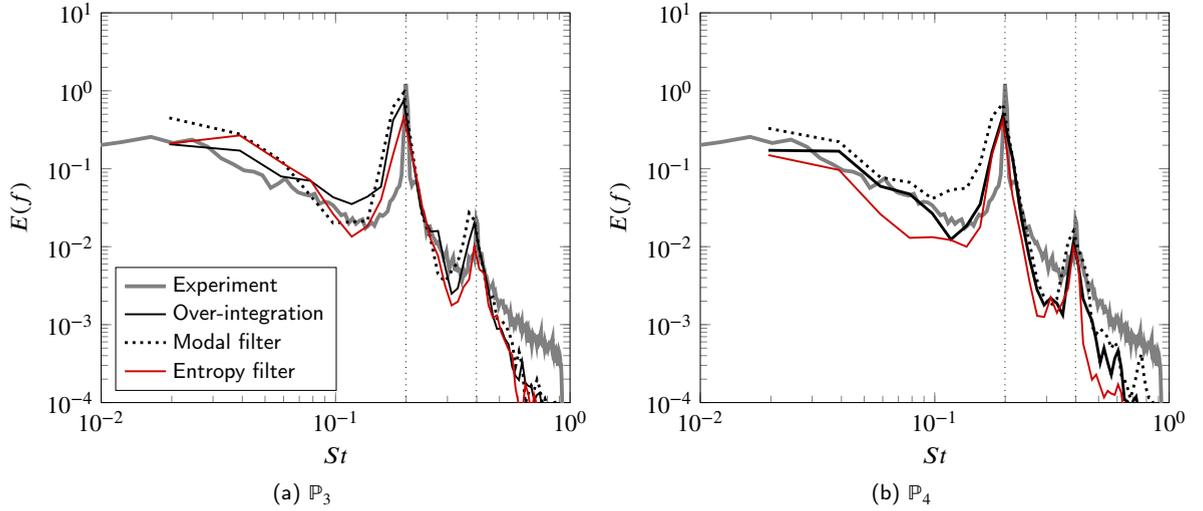

The PSD profiles as computed by the three anti-aliasing approaches are shown in \cref{fig:psd} for both the $\mathbb P_3$ and $\mathbb P_4$ approximations. The experimental results show distinct primary and secondary peaks in the PSD at Strouhal numbers of $St = 0.1994$ and $St = 0.3987$, respectively. For all anti-aliasing approaches at both approximation orders, the primary peak was well predicted in comparison to the experimental results, within the sampling rate error of the periodogram. It must be noted though that at $\mathbb P_3$, the entropy filtering approach showed the most distinct primary peak that was most similar to the experimental results. In comparison, the over-integration and modal filtering approaches showed a broader primary peak. At $\mathbb P_4$, both the entropy filtering and over-integration results showed very similar primary peaks, whereas the primary peak of the modal filtering results was not as distinct. For the secondary peak, some variation in the results was observed. Much like with the primary peak, the over-integration and entropy-filtering results showed very similar predictions in both the location and the prominence of the secondary peak, with good agreement with the experimental results. This observation extended to both the $\mathbb P_3$ and $\mathbb P_4$ approximations. However, the modal filtering approach underpredicted the frequency of the secondary peak at $\mathbb P_3$, but this underprediction was not evident in the $\mathbb P_4$ results.

\subsection{Computational cost}
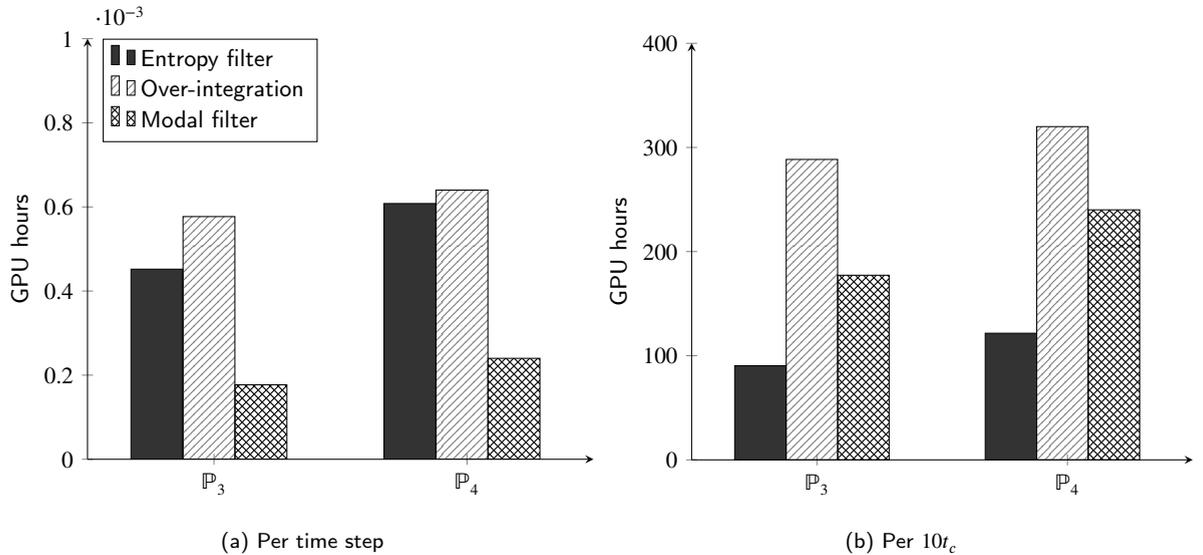
\begin{figure}[htbp!]
    \centering
    \subfloat[Per time step]{\adjustbox{width=0.48\linewidth, valign=b}{\input{figs/cost2}}}
    \subfloat[Per $10 t_c$]{\adjustbox{width=0.48\linewidth, valign=b}{\input{figs/cost1}}}
    \caption{\label{fig:cost} Wall clock time in NVIDIA V100 GPU hours per time step (left) and per 10 flows over chord (right) using a $\mathbb P_3$ and $\mathbb P_4$ approximation with varying anti-aliasing methods. Simulations performed on 16 NVIDIA V100 GPUs. 
    }
\end{figure}

Finally, a comparison of the computational cost of the various approaches was performed with respect to both the wall-clock time required to compute one time step and 10 flows over chord, respectively. The comparison was performed across 16 NVIDIA V100 GPUs, with the results shown in \cref{fig:cost}. As expected due to the large bandwidth and compute requirements of evaluating the projection, the cost of the over-integration approach was the highest per time step, requiring $27\%$ and $225\%$ more compute time than the entropy filtering and modal filtering approaches, respectively, at $\mathbb P_3$. At $\mathbb P_4$, the computational cost of the entropy filtering approach increased proportionally more than the other two approaches, such that the cost of the over-integration approach was only $5\%$ more than the entropy filtering approach per time step but still $160\%$ more than the modal filtering approach. 

However, while the over-integration approach was the most costly and the modal filtering approach was the last costly \textit{per time step}, the total cost of the three approaches varied due to the different time step restrictions of the respective methods. Due to the larger possible time step as a result of the increased nonlinear stability properties of entropy filtering, the total cost of the entropy filtering approach was the lowest at both $\mathbb P_3$ and $\mathbb P_4$, requiring 90.4 GPU hours and 121.6 GPU hours, respectively, to compute 10 flows over chord. In contrast, the cost of the over-integration approach was still the highest, requiring 288.5 GPU hours and 320 GPU hours, respectively, while the cost of the modal filtering approach was approximately in between the two. As such, while the entropy filtering approach may be costly to evaluate per time step, it is possible that it may still decrease the overall computational cost of simulating flows which may exhibit numerical instabilities. 

%% file: figs/naca_p3_cp.tex
    \begin{tikzpicture}[spy using outlines={rectangle, height=3cm,width=2.3cm, magnification=3, connect spies}]
		\begin{axis}[name=plot1,
		    axis line style={latex-latex},
		    axis x line=left,
            axis y line=left,
            clip mode=individual,
		    xlabel={$x$},
		    xtick={0,0.2,0.4,0.6,0.8,1},
    		xmin=-0.02,
    		xmax=1.02,
    	    x tick label style={
        		/pgf/number format/.cd,
            	fixed,
            	fixed zerofill,
            	precision=1,
        	    /tikz/.cd},
    		ylabel={$C_p$},
    		ylabel style={rotate=-90},
    		ytick={-1.5, -1.0, -0.5, 0.0, 0.5, 1.0, 1.5, 2.0},
    		ymin=-1.8,
    		ymax=1.5,
    		y dir=reverse,
            y axis line style = {stealth-},
    		y tick label style={
        		/pgf/number format/.cd,
            	fixed,
            	fixed zerofill,
            	precision=1,
        	    /tikz/.cd},
    		legend style={at={(1.0, 0.03)},anchor=south east,font=\small},
    		legend cell align={left},
    		style={font=\normalsize}]
    		
			\addplot[color=black, style={thin}, only marks, mark=o, mark options={scale=1.1}, mark repeat = 1, mark phase = 0]
				table[x=x,y=cp,col sep=comma,unbounded coords=jump]{./figs/data/naca/ref_exp.csv};
    		\addlegendentry{Exp.}
   
			\addplot[color=black, style={very thick}]
				table[x=x,y=cp,col sep=comma,unbounded coords=jump]{./figs/data/naca/naca0021p3aa_cp.csv};
			\addlegendentry{OI}
			
			\addplot[color=black, style={ultra thick, dotted}]
				table[x=x,y=cp,col sep=comma,unbounded coords=jump]{./figs/data/naca/naca0021p3mf_cp.csv};
			\addlegendentry{MF}

			\addplot[color=red!80!black, style={very thick}]
				table[x=x,y=cp,col sep=comma,unbounded coords=jump]{./figs/data/naca/naca0021p3ef_cp.csv};
			\addlegendentry{EF}
			
		\end{axis} 		
	\end{tikzpicture}

%% file: figs/naca_p4_cp.tex
    \begin{tikzpicture}[spy using outlines={rectangle, height=3cm,width=2.3cm, magnification=3, connect spies}]
		\begin{axis}[name=plot1,
		    axis line style={latex-latex},
		    axis x line=left,
            axis y line=left,
            clip mode=individual,
		    xlabel={$x$},
		    xtick={0,0.2,0.4,0.6,0.8,1},
    		xmin=-0.02,
    		xmax=1.02,
    	    x tick label style={
        		/pgf/number format/.cd,
            	fixed,
            	fixed zerofill,
            	precision=1,
        	    /tikz/.cd},
    		ylabel={$C_p$},
    		ylabel style={rotate=-90},
    		ytick={-1.5, -1.0, -0.5, 0.0, 0.5, 1.0, 1.5, 2.0},
    		ymin=-1.8,
    		ymax=1.5,
    		y dir=reverse,
            y axis line style = {stealth-},
    		y tick label style={
        		/pgf/number format/.cd,
            	fixed,
            	fixed zerofill,
            	precision=1,
        	    /tikz/.cd},
    		legend style={at={(1.0, 0.03)},anchor=south east,font=\small},
    		legend cell align={left},
    		style={font=\normalsize}]
    		
			\addplot[color=black, style={thin}, only marks, mark=o, mark options={scale=1.1}, mark repeat = 1, mark phase = 0]
				table[x=x,y=cp,col sep=comma,unbounded coords=jump]{./figs/data/naca/ref_exp.csv};

			\addplot[color=black, style={very thick}]
				table[x=x,y=cp,col sep=comma,unbounded coords=jump]{./figs/data/naca/naca0021p4aa_cp.csv};
			
			\addplot[color=black, style={ultra thick, dotted}]
				table[x=x,y=cp,col sep=comma,unbounded coords=jump]{./figs/data/naca/naca0021p4mf_cp.csv};
			
			\addplot[color=red!80!black, style={very thick}]
				table[x=x,y=cp,col sep=comma,unbounded coords=jump]{./figs/data/naca/naca0021p4ef_cp.csv}; 
			
		\end{axis} 		
	\end{tikzpicture}

%% file: figs/u_x1_p3.tex
     \begin{tikzpicture}[spy using outlines={rectangle, height=3cm,width=2.3cm, magnification=3, connect spies}]
		\begin{axis}[name=plot1,
		    axis line style={latex-latex},
		    axis x line=left,
            axis y line=left,
            clip mode=individual,
		    xlabel={$y/D$},
		    xtick={-1.5, -1, 0, 1, 2},
    		xmin=-1.5,
    		xmax=1.5,
    		x tick label style={
        		/pgf/number format/.cd,
            	fixed,
            	precision=1,
        	    /tikz/.cd},
    		ylabel={$u/U_\infty$},
		    xtick={-1.5, -1, -0.5, 0, 0.5, 1.0, 1.5},
    		ymin=-0.5,
    		ymax=1.5,
    		y tick label style={
        		/pgf/number format/.cd,
            	fixed,
            	precision=1,
        	    /tikz/.cd},
    		legend style={at={(1.0,0.03)},anchor=south east,font=\small},
    		legend cell align={left},
    		style={font=\normalsize}]

			\addplot[color=black, style={very thick}]
				table[x=y,y=u_x1,col sep=comma,unbounded coords=jump]{./figs/data/naca/p3aa_profiles.csv};
			\addlegendentry{OI}
			\addplot[color=black, style={ultra thick, dotted}]
				table[x=y,y=u_x1,col sep=comma,unbounded coords=jump]{./figs/data/naca/p3mf_profiles.csv};
			\addlegendentry{MF}
			\addplot[color=red!80!black, style={very thick}]
				table[x=y,y=u_x1,col sep=comma,unbounded coords=jump]{./figs/data/naca/p3ef_profiles.csv};
			\addlegendentry{EF}
		\end{axis}
	\end{tikzpicture}

%% file: figs/u_x2_p3.tex
     \begin{tikzpicture}[spy using outlines={rectangle, height=3cm,width=2.3cm, magnification=3, connect spies}]
		\begin{axis}[name=plot1,
		    axis line style={latex-latex},
		    axis x line=left,
            axis y line=left,
            clip mode=individual,
		    xlabel={$y/D$},
		    xtick={-1.5, -1, -0.5, 0, 0.5, 1.0, 1.5},
    		xmin=-1.5,
    		xmax=1.5,
    		x tick label style={
        		/pgf/number format/.cd,
            	fixed,
            	precision=1,
        	    /tikz/.cd},
    		ylabel={$u/U_\infty$},
    		ytick={-0.5, 0, 0.5, 1.0, 1.5},
    		ymin=-0.5,
    		ymax=1.5,
    		y tick label style={
        		/pgf/number format/.cd,
            	fixed,
            	precision=1,
        	    /tikz/.cd},
    		legend style={at={(1.0,0.03)},anchor=south east,font=\small},
    		legend cell align={left},
    		style={font=\normalsize}]

			\addplot[color=black, style={very thick}]
				table[x=y,y=u_x2,col sep=comma,unbounded coords=jump]{./figs/data/naca/p3aa_profiles.csv}; 
   
			\addplot[color=black, style={ultra thick, dotted}]
				table[x=y,y=u_x2,col sep=comma,unbounded coords=jump]{./figs/data/naca/p3mf_profiles.csv};
   
			\addplot[color=red!80!black, style={very thick}]
				table[x=y,y=u_x2,col sep=comma,unbounded coords=jump]{./figs/data/naca/p3ef_profiles.csv};
    		    		
		\end{axis}
	\end{tikzpicture}

%% file: figs/v_x1_p3.tex
     \begin{tikzpicture}[spy using outlines={rectangle, height=3cm,width=2.3cm, magnification=3, connect spies}]
		\begin{axis}[name=plot1,
		    axis line style={latex-latex},
		    axis x line=left,
            axis y line=left,
            clip mode=individual,
		    xlabel={$y/D$},
		    xtick={-1.5, -1, -0.5, 0, 0.5, 1.0, 1.5},
    		xmin=-1.5,
    		xmax=1.5,
    		x tick label style={
        		/pgf/number format/.cd,
            	fixed,
            	precision=1,
        	    /tikz/.cd},
    		ylabel={$v/U_\infty$},
    		ytick={-0.3, -0.2, -0.1, 0, 0.1, 0.2, 0.3},
    		ymin=-0.3,
    		ymax=0.3,
    		y tick label style={
        		/pgf/number format/.cd,
            	fixed,
            	precision=1,
        	    /tikz/.cd},
    		legend style={at={(1.0,0.03)},anchor=south east,font=\small},
    		legend cell align={left},
    		style={font=\normalsize}]

			\addplot[color=black, style={very thick}]
				table[x=y,y=v_x1,col sep=comma,unbounded coords=jump]{./figs/data/naca/p3aa_profiles.csv}; 
			\addplot[color=black, style={ultra thick, dotted}]
				table[x=y,y=v_x1,col sep=comma,unbounded coords=jump]{./figs/data/naca/p3mf_profiles.csv}; 
			\addplot[color=red!80!black, style={very thick}]
				table[x=y,y=v_x1,col sep=comma,unbounded coords=jump]{./figs/data/naca/p3ef_profiles.csv} ;

		\end{axis}

	\end{tikzpicture}

%% file: figs/v_x2_p3.tex
     \begin{tikzpicture}[spy using outlines={rectangle, height=3cm,width=2.3cm, magnification=3, connect spies}]
		\begin{axis}[name=plot1,
		    axis line style={latex-latex},
		    axis x line=left,
            axis y line=left,
            clip mode=individual,
		    xlabel={$y/D$},
		    xtick={-1.5, -1, -0.5, 0, 0.5, 1.0, 1.5},
    		xmin=-1.5,
    		xmax=1.5,
    		x tick label style={
        		/pgf/number format/.cd,
            	fixed,
            	precision=1,
        	    /tikz/.cd},
    		ylabel={$v/U_\infty$},
    		ytick={-0.3, -0.2, -0.1, 0, 0.1, 0.2, 0.3},
    		ymin=-0.3,
    		ymax=0.3,
    		y tick label style={
        		/pgf/number format/.cd,
            	fixed,
            	precision=1,
        	    /tikz/.cd},
    		legend style={at={(1.0,0.03)},anchor=south east,font=\small},
    		legend cell align={left},
    		style={font=\normalsize}]

			\addplot[color=black, style={very thick}]
				table[x=y,y=v_x2,col sep=comma,unbounded coords=jump]{./figs/data/naca/p3aa_profiles.csv}; 
			\addplot[color=black, style={ultra thick, dotted}]
				table[x=y,y=v_x2,col sep=comma,unbounded coords=jump]{./figs/data/naca/p3mf_profiles.csv}; 
			\addplot[color=red!80!black, style={very thick}]
				table[x=y,y=v_x2,col sep=comma,unbounded coords=jump]{./figs/data/naca/p3ef_profiles.csv}; 
    		    		
		\end{axis}
	\end{tikzpicture}

%% file: figs/u_x1_p4.tex
     \begin{tikzpicture}[spy using outlines={rectangle, height=3cm,width=2.3cm, magnification=3, connect spies}]
		\begin{axis}[name=plot1,
		    axis line style={latex-latex},
		    axis x line=left,
            axis y line=left,
            clip mode=individual,
		    xlabel={$y/D$},
		    xtick={-1.5, -1, 0, 1, 2},
    		xmin=-1.5,
    		xmax=1.5,
    		x tick label style={
        		/pgf/number format/.cd,
            	fixed,
            	precision=1,
        	    /tikz/.cd},
    		ylabel={$u/U_\infty$},
		    xtick={-1.5, -1, -0.5, 0, 0.5, 1.0, 1.5},
    		ymin=-0.5,
    		ymax=1.5,
    		y tick label style={
        		/pgf/number format/.cd,
            	fixed,
            	precision=1,
        	    /tikz/.cd},
    		legend style={at={(1.0,0.03)},anchor=south east,font=\small},
    		legend cell align={left},
    		style={font=\normalsize}]

			\addplot[color=black, style={very thick}]
				table[x=y,y=u_x1,col sep=comma,unbounded coords=jump]{./figs/data/naca/p4aa_profiles.csv};
			\addlegendentry{OI}
			\addplot[color=black, style={ultra thick, dotted}]
				table[x=y,y=u_x1,col sep=comma,unbounded coords=jump]{./figs/data/naca/p4mf_profiles.csv};
			\addlegendentry{MF}
			\addplot[color=red!80!black, style={very thick}]
				table[x=y,y=u_x1,col sep=comma,unbounded coords=jump]{./figs/data/naca/p4ef_profiles.csv};
			\addlegendentry{EF}
		\end{axis}
	\end{tikzpicture}

%% file: figs/u_x2_p4.tex
     \begin{tikzpicture}[spy using outlines={rectangle, height=3cm,width=2.3cm, magnification=3, connect spies}]
		\begin{axis}[name=plot1,
		    axis line style={latex-latex},
		    axis x line=left,
            axis y line=left,
            clip mode=individual,
		    xlabel={$y/D$},
		    xtick={-1.5, -1, -0.5, 0, 0.5, 1.0, 1.5},
    		xmin=-1.5,
    		xmax=1.5,
    		x tick label style={
        		/pgf/number format/.cd,
            	fixed,
            	precision=1,
        	    /tikz/.cd},
    		ylabel={$u/U_\infty$},
    		ytick={-0.5, 0, 0.5, 1.0, 1.5},
    		ymin=-0.5,
    		ymax=1.5,
    		y tick label style={
        		/pgf/number format/.cd,
            	fixed,
            	precision=1,
        	    /tikz/.cd},
    		legend style={at={(1.0,0.03)},anchor=south east,font=\small},
    		legend cell align={left},
    		style={font=\normalsize}]

			\addplot[color=black, style={very thick}]
				table[x=y,y=u_x2,col sep=comma,unbounded coords=jump]{./figs/data/naca/p4aa_profiles.csv}; 
   
			\addplot[color=black, style={ultra thick, dotted}]
				table[x=y,y=u_x2,col sep=comma,unbounded coords=jump]{./figs/data/naca/p4mf_profiles.csv};
   
			\addplot[color=red!80!black, style={very thick}]
				table[x=y,y=u_x2,col sep=comma,unbounded coords=jump]{./figs/data/naca/p4ef_profiles.csv};
    		    		
		\end{axis}
	\end{tikzpicture}

%% file: figs/v_x1_p4.tex
     \begin{tikzpicture}[spy using outlines={rectangle, height=3cm,width=2.3cm, magnification=3, connect spies}]
		\begin{axis}[name=plot1,
		    axis line style={latex-latex},
		    axis x line=left,
            axis y line=left,
            clip mode=individual,
		    xlabel={$y/D$},
		    xtick={-1.5, -1, -0.5, 0, 0.5, 1.0, 1.5},
    		xmin=-1.5,
    		xmax=1.5,
    		x tick label style={
        		/pgf/number format/.cd,
            	fixed,
            	precision=1,
        	    /tikz/.cd},
    		ylabel={$v/U_\infty$},
    		ytick={-0.3, -0.2, -0.1, 0, 0.1, 0.2, 0.3},
    		ymin=-0.3,
    		ymax=0.3,
    		y tick label style={
        		/pgf/number format/.cd,
            	fixed,
            	precision=1,
        	    /tikz/.cd},
    		legend style={at={(1.0,0.03)},anchor=south east,font=\small},
    		legend cell align={left},
    		style={font=\normalsize}]

			\addplot[color=black, style={very thick}]
				table[x=y,y=v_x1,col sep=comma,unbounded coords=jump]{./figs/data/naca/p4aa_profiles.csv}; 
			\addplot[color=black, style={ultra thick, dotted}]
				table[x=y,y=v_x1,col sep=comma,unbounded coords=jump]{./figs/data/naca/p4mf_profiles.csv}; 
			\addplot[color=red!80!black, style={very thick}]
				table[x=y,y=v_x1,col sep=comma,unbounded coords=jump]{./figs/data/naca/p4ef_profiles.csv};

		\end{axis}

	\end{tikzpicture}

%% file: figs/v_x2_p4.tex
     \begin{tikzpicture}[spy using outlines={rectangle, height=3cm,width=2.3cm, magnification=3, connect spies}]
		\begin{axis}[name=plot1,
		    axis line style={latex-latex},
		    axis x line=left,
            axis y line=left,
            clip mode=individual,
		    xlabel={$y/D$},
		    xtick={-1.5, -1, -0.5, 0, 0.5, 1.0, 1.5},
    		xmin=-1.5,
    		xmax=1.5,
    		x tick label style={
        		/pgf/number format/.cd,
            	fixed,
            	precision=1,
        	    /tikz/.cd},
    		ylabel={$v/U_\infty$},
    		ytick={-0.3, -0.2, -0.1, 0, 0.1, 0.2, 0.3},
    		ymin=-0.3,
    		ymax=0.3,
    		y tick label style={
        		/pgf/number format/.cd,
            	fixed,
            	precision=1,
        	    /tikz/.cd},
    		legend style={at={(1.0,0.03)},anchor=south east,font=\small},
    		legend cell align={left},
    		style={font=\normalsize}]

			\addplot[color=black, style={very thick}]
				table[x=y,y=v_x2,col sep=comma,unbounded coords=jump]{./figs/data/naca/p4aa_profiles.csv}; 
			\addplot[color=black, style={ultra thick, dotted}]
				table[x=y,y=v_x2,col sep=comma,unbounded coords=jump]{./figs/data/naca/p4mf_profiles.csv}; 
			\addplot[color=red!80!black, style={very thick}]
				table[x=y,y=v_x2,col sep=comma,unbounded coords=jump]{./figs/data/naca/p4ef_profiles.csv}; 
    		    		
		\end{axis}
	\end{tikzpicture}

%% file: figs/naca_p3_psd.tex
\begin{tikzpicture}[spy using outlines={rectangle, height=3cm,width=2.3cm, magnification=3, connect spies}]
		\begin{loglogaxis}[name=plot1,
			xlabel={$St$},
    		xmin=1e-2,xmax=1,
    		ylabel={$E(f)$},
    		ymin=1e-4,ymax=1e1,
    		legend style={at={(0.03,0.03)},anchor=south west,font=\small},
    		legend cell align={left},
    		style={font=\normalsize}]
    		
			\addplot[color=gray, style={ultra thick}]
				table[x=f,y=p,col sep=comma,unbounded coords=jump]{./figs/data/naca/psd_exp.csv};
			\addlegendentry{Experiment}
			
			\addplot[color=black, style={thick}]
				table[x=f,y=p,col sep=comma,unbounded coords=jump]{./figs/data/naca/psd_p3_aa.csv};
			\addlegendentry{Over-integration}
			
			\addplot[color=black, style={very thick, dotted}]
				table[x=f,y=p,col sep=comma,unbounded coords=jump]{./figs/data/naca/psd_p3_mf.csv};
			\addlegendentry{Modal filter}

			\addplot[color=red!80!black, style={thick}]
				table[x=f,y=p,col sep=comma,unbounded coords=jump]{./figs/data/naca/psd_p3_ef.csv};
			\addlegendentry{Entropy filter}

			\addplot[color=black!50!black, style={dotted, thin},forget plot] coordinates{(0.1994, 1e-4) (0.1994, 1e1)};
			
			\addplot[color=black!50!black, style={dotted, thin},forget plot] coordinates{(0.3987, 1e-4) (0.3987, 1e1)};
			
		\end{loglogaxis} 		
	\end{tikzpicture}

%% file: figs/naca_p4_psd.tex
\begin{tikzpicture}[spy using outlines={rectangle, height=3cm,width=2.3cm, magnification=3, connect spies}]
		\begin{loglogaxis}[name=plot1,
			xlabel={$St$},
    		xmin=1e-2,xmax=1,
    		ylabel={$E(f)$},
    		ymin=1e-4,ymax=1e1,
    		legend style={at={(0.03,0.03)},anchor=south west,font=\small},
    		legend cell align={left},
    		style={font=\normalsize}]
    		
			\addplot[color=gray, style={ultra thick}]
				table[x=f,y=p,col sep=comma,unbounded coords=jump]{./figs/data/naca/psd_exp.csv};
			
			\addplot[color=black, style={very thick}]
				table[x=f,y=p,col sep=comma,unbounded coords=jump]{./figs/data/naca/psd_p4_aa.csv};
			
			\addplot[color=black, style={very thick, dotted}]
				table[x=f,y=p,col sep=comma,unbounded coords=jump]{./figs/data/naca/psd_p4_mf.csv};

			\addplot[color=red!80!black, style={thick}]
				table[x=f,y=p,col sep=comma,unbounded coords=jump]{./figs/data/naca/psd_p4_ef.csv};

			\addplot[color=black!50!black, style={dotted, thin},forget plot] coordinates{(0.1994, 1e-4) (0.1994, 1e1)};
			
			\addplot[color=black!50!black, style={dotted, thin},forget plot] coordinates{(0.3987, 1e-4) (0.3987, 1e1)};
			
		\end{loglogaxis} 		
	\end{tikzpicture}

%% file: figs/cost2.tex
\begin{tikzpicture}[spy using outlines={rectangle, height=3cm,width=2.3cm, magnification=3, connect spies}]
\begin{axis} [ybar,
		axis line style={latex-latex},
	    axis x line=left,
        axis y line=left,
        bar width=20pt,
    	xmin=.5, xmax=2.5,
    	ymin=0, ymax=1e-3,
    	ylabel={GPU hours},
        xtick={1, 2},
        xticklabels={$\mathbb P_3$, $\mathbb P_4$},
        clip mode=individual,
    	legend style={at={(0.03, 1.00)},anchor=north west},
    	legend cell align={left}]
     
\addplot[draw=black,fill=black!80]
    coordinates {
    	(1, 4.52e-4) 
    	(2, 6.08e-4) 
    };
\addlegendentry{Entropy filter}

\addplot[draw=black,pattern=north east lines, pattern color = gray]
    coordinates {
    	(0.98, 5.77e-4) 
    	(1.98, 6.4e-4)
    };
\addlegendentry{Over-integration}

\addplot[draw=black,pattern=crosshatch, pattern color = black]
    coordinates {
    	(0.96, 1.77e-4) 
    	(1.96, 2.40e-4)
    };
\addlegendentry{Modal filter}


\end{axis}

\end{tikzpicture}

%% file: figs/cost1.tex
\begin{tikzpicture}[spy using outlines={rectangle, height=3cm,width=2.3cm, magnification=3, connect spies}]
\begin{axis} [ybar,
		axis line style={latex-latex},
	    axis x line=left,
        axis y line=left,
        bar width=20pt,
    	xmin=.5, xmax=2.5,
    	ymin=0, ymax=400,
    	ylabel={GPU hours},
        xtick={1, 2},
        xticklabels={$\mathbb P_3$, $\mathbb P_4$},
        clip mode=individual,
    	legend style={at={(0.03, 1.00)},anchor=north west},
    	legend cell align={left}]
     
\addplot[draw=black,fill=black!80]
    coordinates {
    	(1, 90.4) 
    	(2, 121.6) 
    };

\addplot[draw=black,pattern=north east lines, pattern color = gray]
    coordinates {
    	(0.98, 288.5) 
    	(1.98, 320)
    };

\addplot[draw=black,pattern=crosshatch, pattern color = black]
    coordinates {
    	(0.96, 177.3) 
    	(1.96, 240.1)
    };


\end{axis}

\end{tikzpicture}

%% file: conclusion.tex
\section{Conclusions}\label{sec:conclusion}
In this work, we evaluate the potential of entropy filtering as an anti-aliasing method for high-order discontinuous spectral element approximations of under-resolved turbulent flows. The approach compared to standard over-integration and modal filtering techniques for the under-resolved implict large eddy simulation of a NACA0021 airfoil in deep stall at a Reynolds number of 270,000. It was observed that the entropy filtering approach could robustly mitigate aliasing driven instabilities with accuracy on par with over-integration methods, showing good agreement with experimental data. Additionally, due to the lack of tunable parameters, the results of the entropy filtering approach were more robust and accurate than the modal filtering approach, but showed some degradation in accuracy with increased approximation order. However, the cost of the entropy filtering approach was substantial per time step, on par with the over-integration approach and significantly higher than the modal filtering approach. This additional cost was mitigated by the larger admissible time step possible due to the added nonlinear stability of the entropy filter, such that the total cost of the simulations was lower with the entropy filtering approach in comparison to both over-integration and modal filtering. These results indicate that the entropy filtering method may be an effective anti-aliasing approach for under-resolved turbulent flows. Furthermore, the approach presents a possible unified framework for both shock-capturing and anti-aliasing which may be beneficial for the simulation of high-speed turbulent flows.